\documentclass[UTF8, a4paper, 10pt]{ctexart}
\usepackage[left=2.50cm,  right=2.50cm,  top=2.50cm,  bottom=2.50cm]{geometry} 
\CTEXsetup[format={\Large\bfseries}]{section} 

\usepackage{helvet}

\usepackage[colorlinks,
            linkcolor=blue,
            anchorcolor=blue,
            citecolor=black,
            ]{hyperref}
\usepackage{amsmath,  amsfonts,  amssymb} 
\usepackage[english]{babel}
\usepackage{color}      
\usepackage{graphicx}   
\usepackage{url}        
\usepackage{bm}         
\usepackage{multirow}
\usepackage{booktabs}

\usepackage{epstopdf}
\usepackage{epsfig}
\usepackage{algorithm}
\usepackage{algorithmic}

\usepackage{fancyhdr} 
\usepackage{hyperref} 
\hypersetup{colorlinks,  bookmarks,  unicode} 

\usepackage[marginal]{footmisc}

\usepackage[marginal]{footmisc}

\title {\textbf{An Extension on Neighbor Sum Distinguishing Total Coloring of Graphs}}

\author{\normalsize Jing-zhi \textsc{Chang}$^{1}$,\quad  Chao \textsc{Yang}$^{1,*}$,\quad  Zhi-xiang \textsc{Yin}$^{1}$,\quad  Bing \textsc{Yao}$^{2}$ \\
\normalsize1. School of Mathematics, Physics and Statistics; Center of Intelligent Computing and Applied Statistics,\\
\normalsize Shanghai University of Engineering Science, Shanghai, 201620, China\\
\normalsize2. College of Mathematics and Statistics, Northwest Normal University, Lanzhou, 730070, China}

\date{} 

\begin{document}
    \maketitle
\footnote{\noindent
Supported by National Natural Science Foundation of China under Grant Nos. 61672001, 61662066, 62072296.\\
$^*$Corresponding author. E-mail address: yangchao@sues.edu.cn}

\renewcommand{\abstractname}{} 

\noindent   

\noindent {\textbf{Abstract:}}
Let $f: V(G)\cup E(G)\rightarrow \{1,2,\dots,k\}$ be a non-proper total $k$-coloring of $G$. Define a weight function on total coloring as
$$\phi(x)=f(x)+\sum\limits_{e\ni x}f(e)+\sum\limits_{y\in N(x)}f(y),$$ where $N(x)=\{y\in V(G)|xy\in E(G)\}$.
If $\phi(x)\neq \phi(y)$ for any edge $xy\in E(G)$, then $f$ is
called a neighbor full sum distinguishing total $k$-coloring of
$G$. The smallest value $k$ for which $G$ has such a coloring is
called the neighbor full sum distinguishing total chromatic number
of $G$ and denoted by fgndi$_{\sum}(G)$. The coloring is an extension of neighbor sum distinguishing non-proper total coloring. In this paper we
conjecture that fgndi$_{\sum}(G)\leq 3$ for any connected graph
$G$ of order at least three. We prove that the conjecture is true
for (i) paths and cycles; (ii) 3-regular graphs and (iii) stars, complete graphs, trees,
hypercubes, bipartite graphs and complete $r$-partite graphs. In
particular, complete graphs can achieve the upper bound for the above conjecture.



\section{Introduction \ \ }

All considered graphs are finite, undirected, simple and connected.
Let $[s,t]$ denote the set of nonnegative integers $\{s,s+1,s+2,\ldots,t\}$ and $0\leq s<t$. Let $d_G(v)$
and $\Delta(G)$ (or $\Delta$) denote the degree of vertex $v$ and the maximum degree of $G$, respectively.
Let $d$-vertex denote the vertex of degree $d$, $1\leq d\leq \Delta$. For general theoretic notations, we follow [3].

Graph coloring theory has a wide range of applications in many fields, such as computer science, physics, chemistry and network theory.
Specifically related to time tabling and scheduling, frequency assignment problem, register allocation, computer security, coding theory,
communication network and so on. Since customers have increased dramatically, it yields a confliction between the increasing customers
and the limited expansion of communication network resource. Driven by this background, a class of distinguishing coloring on the sums of colors of vertices and edges has attracted extensive attention.
Karo\'{n}ski et al. [6] firstly introduced and investigated neighbor sum distinguishing edge coloring of graphs, and they proposed a famous 1-2-3 Conjecture. Toward the 1-2-3 Conjecture, Kar\'{o}nski, Luczak and Thomason [7] showed that if $G$ is a $k$-colorable graph
with $k$ odd then $G$ admits a vertex-coloring $k$-edge-weighting. So, for the class of
3-colorable graphs, including bipartite graphs, the answer is affirmative. However,
in general, this question is still open. Addario-Berry et al.[1] showed
that every graph without isolated edges has a proper $k$-weighting when $k=30$. After improvements to $k=15$ in [2] and $k=13$ in [10],
Kalkowski, Kar\'{o}nski, and Pfender [7] showed that every graph without isolated edges has a proper $5$-weighting. Przybylo [9] showed that every $d$-regular graph with
$d\geq 2$ admits a vertex-coloring edge 4-weighting and every $d$-regular graph with $d\ge 10^{8}$ admits a vertex-coloring edge 3-weighting.
Later, Przybylo and Wozniak [8] added the vertex coloring to the weight of edges, they gave the notation of neighbor sum distinguishing total coloring of graphs, meanwhile, they put forward to a 1-2 conjecture with respect to this definition. Thus far it is known that for every graph $G$, tgndi$_{\sum}(G)\leq 3$ (see [5]), where gndi$_{\sum}(G)$ is the neighbor sum distinguishing total chromatic number of $G$.
Recently, Flandrin et al. [4] considered the sum of the colors of neighbors of a vertex based on the neighbor sum distinguishing total coloring, they introduced a new coloring which is called the neighbor full sum distinguishing total coloring, while they didn't give a depth study for this coloring, so we continue to study this type of coloring in this paper.

\vskip 2mm

\noindent {\textbf{Definition 1.}} [4]
\emph{Let $f: V(G)\cup E(G)\rightarrow [1, k]$ be a non-proper $k$-total coloring of $G$.
Set $\phi(x)=f(x)+\sum\limits_{e\ni x}f(e)+\sum\limits_{y\in N(x)}f(y)$, where $N(x)=\{y\in V(G)|xy\in E(G)\}$.
For any edge $xy\in E(G)$, if $\phi(x)\neq \phi(y)$, then $f$ is called a neighbor full sum distinguishing (NFSD) total $k$-coloring of $G$.
The smallest value $k$ for which $G$ has an NFSD-total coloring is called the neighbor full sum distinguishing total chromatic number of $G$ and
denoted by} fgndi$_{\sum}(G)$.

\vskip 2mm

Evidently, when searching for the NFSD-total coloring it is sufficient to
restrict our attention to connected graphs. Observe also, that $G=K_2$ does not have
any NFSD-total coloring. So, we shall consider only connected graphs
with at least three vertices. We propose the following conjecture.

\vskip 2mm

\noindent {\textbf{Conjecture 2.}}\ \ \emph{For every connected graph $G$ and $G$ is not $K_2$,} fgndi$_{\sum}(G)\leq 3$.

\vskip 2mm

By Definition 1, the following result is easy to obtain.

\vskip 2mm

\noindent {\textbf{Lemma 3.}}\ \ \emph{Let $G$ be a connected simple graph of order at least three. Then $(i)$} fgndi$_{\sum}(G)=1$ \emph{if $G$ contains no adjacent $d$-vertices;
and} $(ii)$ fgndi$_{\sum}(G)\geq 2$ \emph{if $G$ contains adjacent $d$-vertices.}

\vskip 2mm

\textbf{\emph{Proof}}
For any two adjacent vertices $u$ and $v$ of $G$, (i) if $G$ contains no adjacent $d$-vertices, namely $d_G(u)\neq d_G(v)$, then we color all vertices and edges of $G$ with 1, and it gets that $\phi(u)=2d_{G}(u)+1\neq 2d_{G}(v)+1=\phi(v)$; (ii) if $d_G(u)=d_G(v)=d$, then fgndi$_{\sum}(G)\geq 2$, otherwise, $u$ and $v$ receive the same weight, a contradiction. $\Box$

\vskip 2mm

We organize the paper as follows. In Section 2, the neighbor full sum distinguishing total chromatic number of paths and cycles are determined.
In Section 3, we offer an important structural lemma that every connected graph $G$ contains
a $m$-partite spanning subgraph $H$ such that
$(1-\frac{1}{m})d_G(v)\leq d_H(v)$. Therefore, every 3-regular
graph $G$ has a maximal bipartite spanning subgraph $H$ such that
$G-E(H)$ is either isolated vertices or isolated edges. Via the
structural between $H$ and $G-E(H)$ of 3-regular graphs $G$ and combining with a coloring
algorithm, we get that fgndi$_{\Sigma}(G)\leq 3$ for any 3-regular graph $G$. In Section 4, we obtain the parameter fgndi$_{\Sigma}(G)$ of several types of graphs with maximum degree $\Delta\geq 4$.

\section{Graphs with $\Delta=2$ \ \ }

\noindent {\textbf{Proposition 4.}}\ \ \emph{Let $P_n$ be a path of order $n$ \emph{($\geq 3$)}.
Then} fgndi$_{\sum}(P_n)=2$ \emph{if} $n\geq 4$ \emph{and} fgndi$_{\sum}(P_3)=1$.

\textbf{\emph{Proof}}
Let~$P_n=x_1x_2\dots x_n$. It is easy to verify that fgndi$_{\Sigma}(P_3)=1$.
By Lemma 3, fgndi$_{\Sigma}(P_n)\geq 2$ for~$n\geq 4$. We define a total coloring $f$: $V(G) \cup E(G)\rightarrow \{ 1,2\}$ as follows:
$$
f\left( x_i \right) =\left\{ \begin{array}{l}
    1~~~if\ i\equiv 1\left(mod~2 \right),\\
    2~~~if\ i\equiv 0\left(mod~2 \right) .
\end{array}\ \right.
$$
And all edges of $P_n$ are colored by $1$.

Taking advantage of the above coloring $f$, we
have~$\phi(x_1)=\phi(x_n)=4$. For any vertex $x_k$ ($2\leq k\leq
n-1$), $\phi(x_k)=6$ if $k$ is even and $\phi(x_k)=7$ if $k$ is
odd, which deduces that $f$ is an NFSD-$2$-total coloring of $P_n$. $\Box$

\vskip 2mm

\noindent {\textbf{Proposition 5.}} \emph{Let $C_n$ be a cycle with order $n (\geq 3)$. Then}
$$
\ fgndi_{\Sigma}(C_n)=\left\{ \begin{array}{l}
    3~~~if\ n=3,\\
    2~~~if\ n\geq 4.
\end{array} \right.
$$

\textbf{\emph{Proof}}  Let~$C_n = x_1x_2\dots x_nx_1$. Clearly, fgndi$_{\Sigma}(C_3) = 3$.
By Lemma 3, fgndi$_{\Sigma}(C_n)\geq 2$ for $n\geq 4$. The following two cases imply that $C_n$ has an NFSD-$2$-total coloring.

\textbf{Case~1.} $n \equiv 1(mod~2)$.

We define a total $2$-coloring $f$ of $C_n$ as below.

$f(x_1)=2$;

$\ f(x_i) =\left\{ \begin{array}{l}
    1,\ i\equiv 1( mod~2) \\
    2,\ i\equiv 0( mod~2 )
\end{array} \right.,\ i\in [2,n]$;

$f(x_1x_2)=f(x_nx_1)=2$;

$f(x_jx_{j+1})=1, j\in [2,n-1]$.

Then
$\phi(x_1) =9$, $\phi(x_2) =8$, $\phi(x_n) = 8$, $\phi(x_{i}) = 7$ if $i$ is odd and $i\in[3, n-1]$,
$\phi(x_{i}) = 6$ if $i$ is even and $i\in[3, n-1]$.
Therefore, $\phi(u) \neq \phi(v)$ for any edge $uv\in E(C_n)$, namely, $f$ is an
NFSD-total $2$-coloring of $C_n$, and thus fgndi$_{\Sigma}(C_n) = 2$.

\textbf{Case~2.} $n \equiv 0(mod~2)$.

We define a total $2$-coloring $f$ of $C_n$ as below.
$$
f\left( x_i \right) =\left\{ \begin{array}{l}
    1,\ i\equiv 1\left( mod~2 \right) \\
    2,\ i\equiv 0\left( mod~2 \right)
\end{array} \right. ,\ i\in [1,n];
$$
Meanwhile, all edges of $C_n$ are colored by $1$.
Then~$\phi(x_{2k})=6\neq \phi(x_{2k+1})=7$. Therefore, $f$ is an
NFSD-total $2$-coloring of $C_n$, and hence fgndi$_{\Sigma}(C_n) =
2$. $\Box$

\section{3-regular graphs\ \ }

This section we investigate fgndi$_{\Sigma}(G)$ of 3-regular graphs. The following lemma is
very crucial to the proof of the main theorem.

\vskip 2mm

\noindent {\textbf{Lemma 6.}}\ \ \emph{Let $G$ be a graph on $n$ vertices. Then it exists a $m$-partite spanning subgraph $H$ such that
$(1-\frac{1}{m})d_G(v)\leq d_H(v)$ for all $v\in V(G)$, where $m$ is a positive
integer and $m\leq n$.}

\vskip 2mm

\textbf{\emph{Proof}}
Let $H$ be a maximal $m$-partite spanning subgraph of $G$ with the greatest
possible number of edges. Let $\{V_1,V_2,\dots,V_m\}$ be the $m$-partition of $V(H)$ and let
$v\in V_1$, $d_{V_i}(v)=|N_{V_i}(v)|$, $N_{V_i}(v)=\{u:u\in V_i,uv\in E(G)\}$,
$i=1,2,\dots,m$. Then $d_{V_1}(v)\leq d_{V_i}(v)$, $i=1,2,\dots,m$. Otherwise, it exists an $i_0$ such that
$d_{V_1}(v)>d_{V_{i_0}}(v)$, and we use $V_1\setminus \{v\}$, $V_{i_0}\cup \{v\}$ instead of $V_1$,
$V_{i_0}$, respectively, and then it generates a new maximal $m$-partite spanning subgraph $H^{'}$ of
$G$. Obviously, $\varepsilon (H^{'})> \varepsilon (H)$, a contradiction. Therefore,

$$(m-1)d_{V_1}(v)\leq \sum_{i=2}^{m}d_{V_i}(v)=d_H(v), $$
~where~$d_H(v)=|N_H(v)|$.

Let $d_G(v)=|N_G(v)|$. Then

$$d_G(v)=d_{V_1}(v)+d_H(v)\leq \frac{1}{m-1}d_H(v)+d_H(v)=\frac{m}{m-1}d_H(v)$$

Hence
$$(1-\frac{1}{m})d_G(v)\leq d_H(v)$$ for all $v\in V(G)$. $\Box$
 \vskip 2mm

Lemma 6 implies that every 3-regular graph $G$ contains a maximal
bipartite spanning subgraph $H$ such that $G-E(H)$ is either
isolated vertices or isolated edges. Via the structural between
$H$ and $G-E(H)$, we further study the neighbor full sum
distinguishing total chromatic number of $3$-regular graphs.

For
the sake of narrative, we fix some natation. Let
$G=(V_{X},V_{Y},E_{X},E_{Y},E_{H})$ be a 3-regular graph with
vertex partition $(V_{X},V_{Y})$ and edge partition
$(E_{X},E_{Y},E_{H})$, where $E_X, E_Y$ and $E_H$ represent the
edge sets in $V_X, V_Y$ and $H$, respectively. The maximal
spanning bipartite subgraph $H$ is the graph with vertex set
$\{x_{i}:1\le i\le m\}\cup \{y_{i} : 1\le i\le n\}$ and edge set
$\{x_{i}y_{i}:1\le i\le m,1\le i\le n\}$. We use $a_{1}$ and
$b_{1}$ to denote the number of vertices with degree 2 and 3 in
$X$, respectively. Let $a_{2}$ and $b_{2}$ denote the number of
vertices with degree 2 and 3 in $Y$, respectively. We use
$e_{x_{i}}$ for $i=1,2,\dots ,\frac{a_{1}}{2}$ to represent the
edge with two endpoints $v_{x_{i}}$, $v'_{x_{i}}$ in $X$. Let
$e_{y_{j}}$ for $j=1,2,\dots ,\frac{a_{2}}{2}$ to represent the
edge with two endpoints $v_{y_{j}}$, $v'_{y_{j}}$ in $Y$.

\vskip 2mm

\noindent {\textbf{Theorem 7.}}\ \ \emph{For any $3$-regular graph $G$,} $2\leq$ fgndi$_{\Sigma}(G)\leq 3$.

\vskip 2mm

\textbf{\emph{Proof}} Let $G$ be a $3$-regular graph. Then fgndi$_{\Sigma}(G)\geq 2$ by Lemma 3.

\textbf{Case 1.}  $a_{1}=a_{2}=0$.

This case implies that $G$ is a 3-regular complete bipartite graph. We color all vertices in $X$ with $1$, color all vertices in $Y$ with 2 and color all edges in $E_{H}$ with 1. Then for any two vertices $v_{x}\in V_{X}$ and $v_{y}\in V_{Y}$, $\phi(v_{x})=10$ and $\phi(v_{y})=8$.

\textbf{Case 2.}  $a_{1}=b_{2}=0$ or $b_{1}=a_{2}=0$.

Without loss of generality, assume that $a_{1}=b_{2}=0$.  We color
all vertices in $X$ with $1$, color all vertices in $Y$ with 2,
color all edges in $E_{H}$ with 1 and color all edges in $E_{Y}$
with 1. Then for any two vertices $v_{x}\in V_{X}$ and $v_{y}\in
V_{Y}$, $\phi(v_{x})=10$ and $\phi(v_{y})=9$. To assure that
$\phi(v_{y_{j}})\ne \phi(v'_{y_{j}})$, recolor an incident edge of
$v_{y_{j}}$ (or $v'_{y_{j}}$) with 3. Then $\phi(v_{x})$ belongs
to $\{10,12,14,16\}$ and $\phi(v_{y})$ is equal to 9 or 11.

\textbf{Case 3.}  $b_{1}=b_{2}=0$.

By Lemma 6, $G$ contains a maximal bipartite spanning subgraph $H$
such that $G-E(H)$ is either isolated vertices or isolated edges.
We color all vertices in $X$ with $1$, color all vertices in $Y$
with 2, color all edges in $E_{H}$ and $E_{X}$ with 1, and color
all edges in $E_{Y}$ with 2. For any edge $v_{x_{i}}v'_{x_{i}}$ in
$X$, select one edge $e_z$ from $E_H$ such that $v_{x_{i}}$ (or
$v'_{x_{i}}$) is an endpoints of $e_z$ and recolor edge $e_z$ with
3, meanwhile, all incident edges (except for $e_z,
v_{x_{i}}v'_{x_{i}}, v_{y_{j}}v'_{y_{j}}$) of
$v_{x_i},v'_{x_i},v_{y_j},v'_{y_j}$ keep the color 1 as before,
and we call these edges being dominated, see Fig.1. Without loss
of generality, assume that $v'_{x_{i}}$ and $v_{y_{j}}$ are
connected by $e_{z}$. Then $\phi(v_{x_{i}})=11,\
\phi(v'_{x_{i}})=9$, $\phi(v_{y_{j}})=12,\ \phi(v'_{y_{j}})=10$.
Continue this procedure $\frac{a_{1}}{2}$ times until the weights
of all adjacent vertices in $G$ are distinct. Now we prove it
feasibility, namely, it verifies that there exists at least one
edge in $E_{H}$ can not be dominated after $\frac{a_{1}}{2}-1$
operations. Suppose that all edges in $H$ are dominated after
$\frac{a_{1}}{2}-1$ operations, and if there still exists a pair
of adjacent vertices $v_{x_{k}}$ and $v'_{x_{k}}$ (or $v_{y_{k}}$
and $v'_{y_{k}}$) having the same weights, then the four incident
edges of $v_{x_{k}}$ and $v'_{x_{k}}$ in $H$ receive the same
color 1. By our coloring rule, it is impossible.

\vskip 2mm

\begin{figure}[h]
\centering
\includegraphics[height=3.5cm]{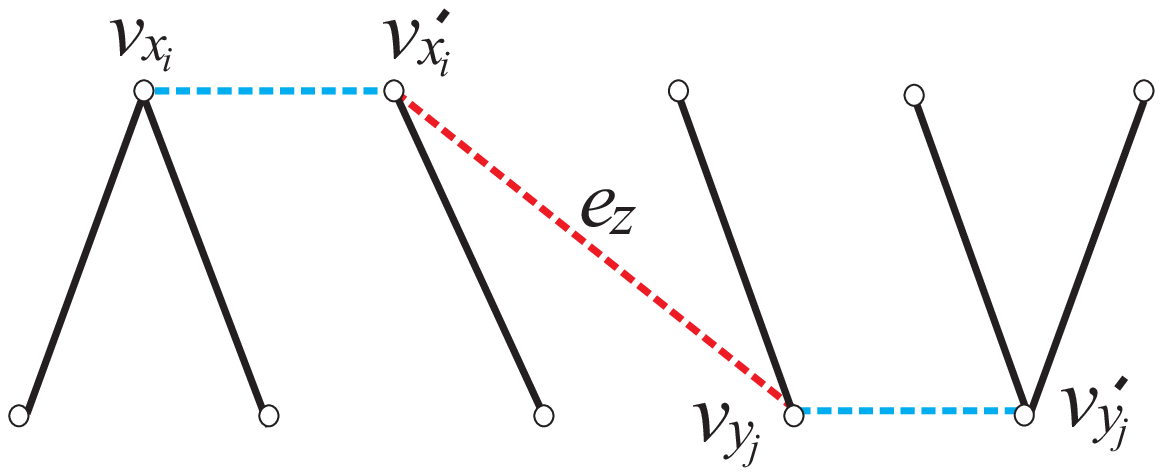}
\begin{center}
{\small Fig.1: Edges labelled by solid line are dominated.}
\end{center}
\end{figure}

\textbf{Case 4.}  $a_{1}\ne a_{2}$, $b_{1}\ne b_{2}$ and they are all positive integers.

We color all vertices in $X$ with $1$, color all vertices in $Y$ with 3, color all edges in $E_{H}$ with 1, color all edges in $E_{X}$ with 2 and color all edges in $E_{Y}$ with 3. In the bipartite graph $H$, let $v_{2}(x)$ and $v_{3}(x)$ be the vertex with degree 2 and 3 in $X$, respectively. Similarly, $v_{2}(y)$ and $v_{3}(y)$ denote the vertex with degree 2 and 3 in $Y$, respectively. The edge between $v_{2}(x)$ and $v_{2}(y)$ is denoted by $\widetilde{e_{2}}$, and the edge connects $v_{2}(x)$ (or $v_{3}(x)$) and $v_{3}(y)$ (or $v_{2}(y)$) is denoted by $\widetilde{e_{2-3}}$.

 Using the technique in Case 3 to change the weight of the vertex which is connected only by $\widetilde{e_{2}}$,
 we can distinguish all adjacent vertices which are connected by $\widetilde{e_{2}}$. But we still need to distinguish adjacent
 vertices which are joined by $\widetilde{e_{2-3}}$. Select an incident edge from $\widetilde{e_{2-3}}$ and color it with 3,
 then it deduces that $\phi(v_{x_{i}})\ne \phi(v'_{x_{i}})$ and $\phi(v_{y_{j}})\ne \phi(v'_{y_{j}})$, meanwhile
 $\phi(v_{2}(x))\in\{12, 14\}$, $\phi(v_{2}(y))\in \{13,15\}$, $\phi(v_{3}(x))\in \{13,15,17,19\}$ and
 $\phi(v_{3}(y))\in \{9,11,13,15\}$. Possibly, there are some cases that the weight of adjacent vertices can not distinguish. We deal with it as follows.

\textbf{Case 4.1.} $\phi(v_{3}(x))=13$. Let
$N(v_{3}(x))=\{v_{y_{1}},\ v_{y_{2}},\ v_{y_{3}}\}$.\par
\textbf{Case 4.1.1.} All vertices in $N(v_{3}(x))$ have the same
weight 13. Recolor $v_{3}(x)$ and its incident edges with 3. Then
$\phi(v_{3}(x))=21$ and
$\phi(v_{y_{1}})=\phi(v_{y_{2}})=\phi(v_{y_{3}})=17$. If one of
$v_{y_{1}},\ v_{y_{2}},\ v_{y_{3}}$ has an adjacent vertex
$v_{x_{l}}$ with weight 17 in $X$, say $v_{y_{1}}$, then recolor
edges $v_{3}(x)v_{y_{1}}$ and  $v_{y_{1}}v_{x_{l}}$ with 2, and we
have $\phi(v_{3}(x))=20$, $\phi(v_{y_{1}})=17$ and
$\phi(v_{x_{l}})=18$.

\textbf{Case 4.1.2.} One of a vertex in $N(v_{3}(x))$ has weight 15, say $v_{y_{1}}$. Recolor edge $v_{3}(x)v_{y_{1}}$ with 2. Then $\phi(v_{3}(x))=14$ and $\phi(v_{y_{1}})=16$.

\textbf{Case 4.1.3.} One of a vertex in $N(v_{3}(x))$ has weight 9, say $v_{y_{2}}$. Recolor edge $v_{3}(x)v_{y_{2}}$ with 2. Then $\phi(v_{3}(x))=14$ and $\phi(v_{y_{2}})=10$.

\textbf{Case 4.1.4.} One of a vertex in $N(v_{3}(x))$ has
weight 11, say $v_{y_{3}}$. Recolor edge $v_{3}(x)v_{y_{3}}$ with
3. Then $\phi(v_{3}(x))=15$ and $\phi(v_{y_{3}})=13$. If there is
a neighbor vertex (say $v_{x_{k}}$) of $v_{y_{3}}$ having weight
13 and a neighbor vertex of $v_{x_{k}}$ having weight 11,
then recolor $v_{3}(x)v_{y_{3}}$ and $v_{y_{3}}v_{x_{k}}$ with 2, and we
have $\phi(v_{3}(x))=\phi(v_{x_{k}})=14$ and $\phi(v_{y_{3}})=13$.
If Case 4.1.4 and Case 4.1.2 appear at the same time, use the
method of Case 4.1.2.

\textbf{Case 4.2.} $\phi(v_{3}(x))=15$.
Let $N(v_{3}(x))=\{v'_{y_{1}},\ v'_{y_{2}},\ v'_{y_{3}}\}$. Then $H$ contains a 2-vertex (vertex of degree 2) with weight 15 in $N(v_{3}(x))$, say $v'_{y_{1}}$.

\textbf{Case 4.2.1.} If $v'_{y_{2}}$ and $v'_{y_{3}}$ have the same weight 15, recolor $v_{3}(x)v'_{y_{2}}$ and $v_{3}(x)v'_{y_{3}}$ with 2, then $\phi(v_{3}(x))=17$ and $\phi(v'_{y_{2}})=\phi(v'_{y_{3}})=16$.

\textbf{Case 4.2.2.} If one of $\{v'_{y_{2}},\ v'_{y_{3}}\}$ has weight 9, assume that $\phi(v'_{y_{2}})=9$.  Recolor edge $v_{3}(x)v'_{y_{2}}$ with 2, then $\phi(v_{3}(x))=16$ and $\phi(v'_{y_{2}})=10$.

\textbf{Case 4.2.3.} If one of $\{v'_{y_{2}},\ v'_{y_{3}}\}$ has weight 11, assume that $\phi(v'_{y_{2}})=11$.  Recolor edge $v_{3}(x)v'_{y_{2}}$ with 3, then $\phi(v_{3}(x))=17$ and $\phi(v'_{y_{2}})=13$.

\textbf{Case 4.2.4.} If one of $\{v'_{y_{2}},\ v'_{y_{3}}\}$ has weight 13, assume that $\phi(v'_{y_{2}})=13$, then two cases appear as follows:

 (i) $d_{H}(v'_{y_{2}})=3$. Recolor edge $v_{3}(x)v'_{y_{2}}$ with 3, then $\phi(v_{3}(x))=17$ and $\phi(v'_{y_{2}})=15$.

 (ii) $d_{H}(v'_{y_{2}})=2$ and $\phi(v'_{y_{3}})=15$. Recolor edge $v_{3}(x)v'_{y_{2}}$ and vertex $v_{3}(x)$ with 3, recolor $v_{3}(x)v'_{y_{1}}$ and $v_{3}(x)v'_{y_{3}}$ with 2, then $\phi(v_{3}(x))=20$ and $\phi(v'_{y_{1}})=16$, $\phi(v'_{y_{2}})=17$, $\phi(v'_{y_{3}})=18$. If $v'_{y_{3}}$ is adjacent to a vertex having weight 18 in $X$, recolor $v_{3}(x)v'_{y_{3}}$ with 1, then $\phi(v_{3}(x))=19$ and $\phi(v'_{y_{3}})=17$. If $v'_{y_{3}}$ is adjacent to a vertex having weight $17$ in $Y$, by our coloring rule, a vertex $v_{0}$ with weight $17$ in $Y$ must have an adjacent vertex $v'_{0}$ with its weight not equal to $18$ in $X$. Then recolor edge $v_{0}v'_{0}$ with 2. It deduces that $\phi(v_{0})=16$ and $\phi(v'_{0})$ reduces 1 than before.

\textbf{Case 4.2.5.} If $v'_{y_{2}}$ and $v'_{y_{3}}$ have the same weight
13 and $d_H(v'_{y_{2}})=d_H(v'_{y_{3}})=2$, we recolor edges
$v_{3}(x)v'_{y_{2}}$ and $v_{3}(x)v'_{y_{3}}$ with 3, recolor edge
$v_{3}(x)v'_{y_{1}}$ and vertex $v_{3}(x)$ with 2 and 3,
respectively, then $\phi(v_{3}(x))=20$, $\phi(v'_{y_{1}})=16$ and
$ \phi(v'_{y_{2}})=\phi(v'_{y_{3}})=17$. If one of $\{v'_{y_{2}},\
v'_{y_{3}}\}$ (say $v'_{y_{2}}$) is adjacent to a vertex
$v'_{x_{k}}$ having weight $17$ in $X$, then recolor edges
$v_{3}(x)v'_{y_{2}}$ and  $v'_{y_{2}}v'_{x_{k}}$ with 2, and it
follows that $\phi(v_{3}(x))=19$, $\phi(v'_{y_{2}})=17$ and
$\phi(v'_{x_{k}})=18$. $\Box$

 \section{Several types of graphs with $\Delta\geq 4$\ \ }

\noindent {\textbf{Proposition 8.}}\ \ \emph{Let $S_n=K_{1,n-1}$ be a star of order $n$.
Then} fgndi$_{\sum}(S_n)=1$.

\textbf{\emph{Proof}}
This conclusion is easily proved by using 1 to color all vertices
and edges of $S_n$. $\Box$

\vskip 2mm

\noindent {\textbf{Proposition 9.}}\ \ \emph{For any complete bipartite graph $K_{m,n}$,
$fgndi_{\sum}(K_{m,n})=1$ if $m\neq n$ and} fgndi$_{\sum}(K_{m,n})=2$ \emph{if $m=n$}.

\textbf{\emph{Proof}}
Suppose that $K_{m,n}=(X, Y, E)$ is a complete bipartite graph
with bipartition classes $X$ and $Y$. Let $|X|=m$ and $|Y|=n$. If
$m\neq n$, then use 1 to color all vertices and edges of
$K_{m,n}$, and it follows that $\phi(x)=2n+1\neq 2m+1=\phi(y)$,
where $x\in X, y\in Y$. For $m=n$, if we use 1 to color all
vertices and edges of $K_{m,n}$, then $\phi(x)=2n+1=2m+1=\phi(y)$
for any edge $xy\in K_{m,n}$, a contradiction, which deduces that
fgndi$_{\sum}(K_{m,n})\geq 2$. We define a total 2-coloring $f$ of
$K_{n,n}$ as follows: using 2 to color each vertex of $Y$ and the
remaining vertices and edges are colored by 1. Then we have
$\phi(x)=3n+1\neq 2n+2=\phi(y)$ for any edge $xy\in K_{n,n}$.
Namely, $f$ is an NFSD-total $2$-coloring of $K_{m,n}$.
$\Box$

\vskip 2mm
\noindent {\textbf{Theorem 10.}}\ \ \emph{For any complete graph $K_n$ $(n\geq 3)$,} fgndi$_{\sum}(K_n)=3$.

\textbf{\emph{Proof}}
It is well known that all vertices are neighbors in $K_n$, so the neighbor full sum distinguishing total coloring is
actually a neighbor sum distinguishing edge coloring of $K_n$. We need only
to consider a neighbor sum distinguishing edge coloring of $K_n$. Suppose that $f$ is a neighbor
sum distinguishing edge 2-coloring of $K_n$ and all vertices of $K_n$ are colored by $1$. For each vertex of $K_n$,
its ($\Delta-1$) incident edges are colored by $1$ and $2$, there exists two vertices $u$ and $v$
such that all incident edges of $u$ are colored by $1$ and all incident edges of $u$ are colored by $2$, a contradiction.
Therefore, fgndi$_{\sum}(K_n)\geq 3$. We offer a method to give a neighbor full sum distinguishing total 3-coloring of $K_n$.

Let $f_3$ be the total coloring of $K_3$ defined as follows:
$f_3(x_i)=1$ for $x_i\in V(K_3), i\in\{1,2,3\}$, $f_3(x_1x_2) = 1, f_3(x_2x_3) = 2, f_3(x_1x_3) = 3$.
Then $\phi(x_1)=3, \phi(x_2)=4, \phi(x_3)=5$.
The coloring $f_n$ will be defined recursively as follows.
If $n$ is odd use 1 to color the vertex $x_{n+1}$, and $1$ to color all edges incident to $x_{n+1}$.
For $n$ even use 1 to color the vertex $x_{n+1}$ and $3$ to color all edges incident to this vertex.

Observe in the coloring $f_n$, if $n$ is odd, then the weights
$\phi(x_i)$ for $i\in [1,n]$ increase by $2$ (with respect to the
weights for $f_{n-1}$) and $\phi(x_{n+1})$ is equal to $2n+1$. If
$n$ is even, then the weights $\phi(x_i)$ for $i\in [1,n]$
increase by $4$ (with respect to the weights for $f_{n-1}$) and
$\phi(x_{n+1})$ is equal to $4n+1$. It follows that $f_n$ is an
NFSD-total $3$-coloring of $K_n$. $\Box$

\vskip 2mm

Theorem 10 implies that there is a type of graphs such that their
NFSD-total chromatic numbers achieving the upper bound of
Conjecture 1. Kalkowski et al. [3] showed that for every connected graph $G$ on order at least three, there exists a coloring of the edges
of $G$ with the colors of [1,5] such that the resulting vertex weighting is a proper vertex coloring of $G$.
Let $K_{n_1,n_2,\dots,n_r}$ be a complete $r$-partite graph with $r$ vertex sets $X_i$ ($i\in [1,r]$) and $|X_i|=n_i$, $\sum_{i=1}^{r}=n$.
If $n_1=n_2=\dots=n_r=\frac{n}{r}$, then
color the edges of $K_{n_1,n_2,\dots,n_r}$ with 5 colors in such a
way that the obtained vertex coloring is proper. Afterwards put 1
on all vertices of $K_{n_1,n_2,\dots,n_r}$, then all weights will
increase by a constant $n-\frac{n}{r}+1$, namely, there exist an
NFSD-total $5$-coloring of $K_{n_1,n_2,\dots,n_r}$. However, the bound 5 can be improved to 3 as follows:

\vskip 2mm

\noindent {\textbf{Theorem 11.}} \emph{Let $K_{n_1,n_2,\dots,n_r}$ be a complete $r$-partite graph with $r$ vertex sets $X_i$ ($i\in [1,r]$) and $|X_i|=n_i$, $\sum_{i=1}^{r}=n$. Then} (i) fgndi$_{\sum}(K_{n_1,n_2,\dots,n_r})=1$ \emph{if} $n_1<n_2<\dots<n_r$; \emph{and} (ii) fgndi$_{\sum}(K_{n_1,n_2,\dots,n_r})\leq 3$ \emph{if} $n_1=n_2=\dots=n_r=\frac{n}{r}$.

\textbf{\emph{Proof}}
Let $K_{n_1,n_2,\dots,n_r}$ be a complete $r$-partite graph with
$r$ vertex sets $X_i$ ($i\in [1,r]$) and $|X_i|=n_i$,
$\sum_{i=1}^{r}=n$.

(i) For $n_1<n_2<\dots<n_r$, we color all vertices
and edges by 1, then $\phi(x_i)=2(n-n_i)+1$, thus $\phi(x_i)\neq
\phi(x_j)$ for $i\neq j$.

(ii) For $n_1=n_2=\dots=n_r=\frac{n}{r}$, we treat each partition $X_i$ of $K_{n_1,n_2,\dots,n_r}$ as a vertex.
Then $K_{n_1,n_2,\dots,n_r}$ can be degenerated into a complete graph $K_r$. Theorem 10 implies that $K_r$ can
achieve an NFSD-total 3-coloring. Let $v_i$ be the vertex of $K_r$. To show that $K_{n_1,n_2,\dots,n_r}$ has an NFSD-total 3-coloring,
we take a method as follows: (1) Vertices of $K_{n_1,n_2,\dots,n_r}$ are colored by 1; (2) Edges between $X_i$ and $X_j$ have the same color with
the edge $v_iv_j$ in $K_r$. Because of the weights of vertices of $K_r$ are different, assume that vertex $v_i$ has weight $f_i$ and $f_1<f_2<\dots<f_r $,
then all vertices in the same part $X_i$ have the same weight $(\frac{n}{r})\cdot f_i+r$
and $(\frac{n}{r})\cdot f_1<(\frac{n}{r})\cdot f_2<\dots<(\frac{n}{r})\cdot f_r$. Therefore, vertices in $r$ partitions have different weights.
$\Box$

\vskip 2mm

\begin{figure}[h]
\centering
\includegraphics[height=5cm]{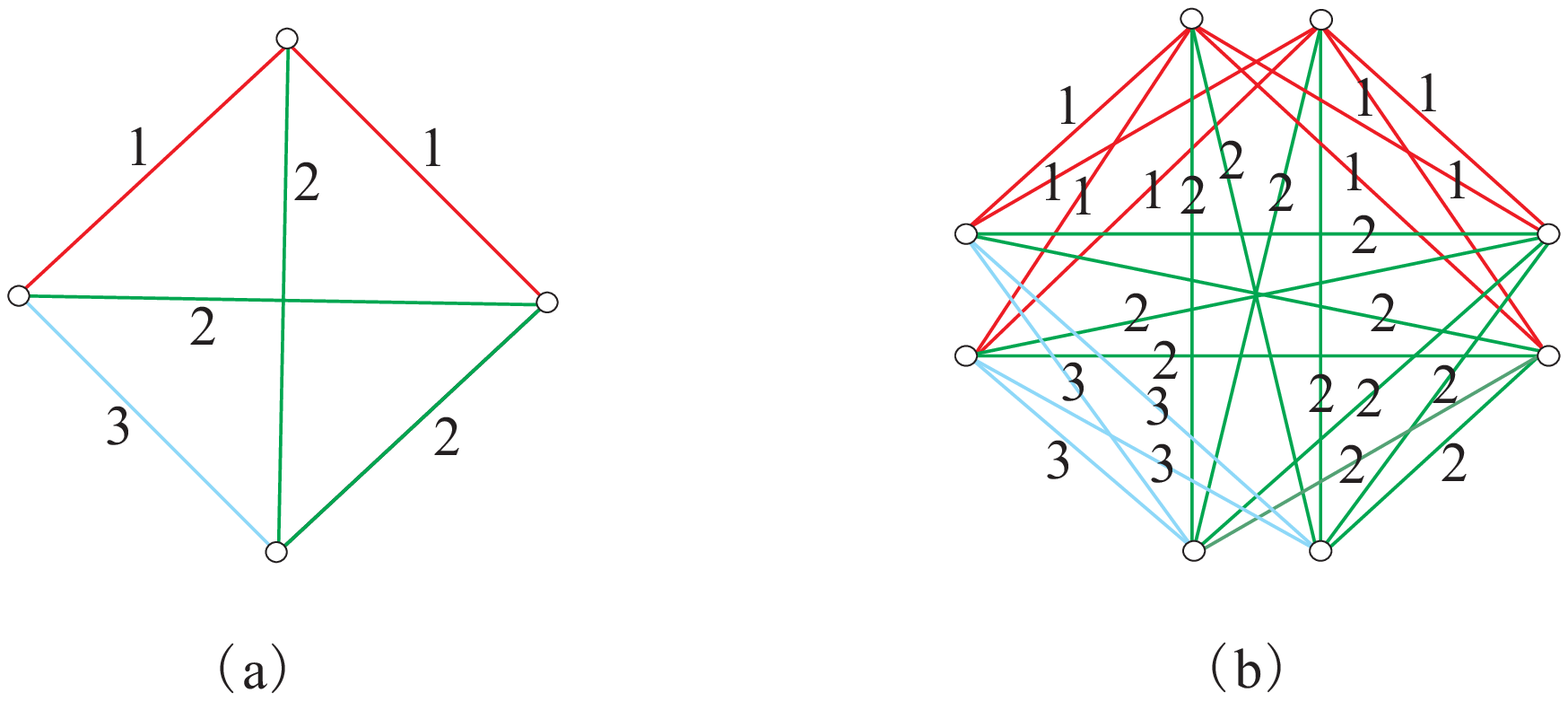}
\begin{center}
{\small Fig.2: (a) $K_4$ and (b) $K_{2,2,2,2}$.}
\end{center}
\end{figure}

\noindent {\textbf{Theorem 12.}}\ \ \emph{For any tree $T$,} (i) fgndi$_{\sum}(T)=1$ \emph{if $T$ contains no adjacent $d$-vertices};
\emph{and} (ii) fgndi$_{\sum}(T)=2$ \emph{if $T$ contains adjacent $d$-vertices}.

\textbf{\emph{Proof}}
By Lemma 3, conclusion (i) is obvious. Next we consider the case that $T$ has adjacent $d$-vertices.
The proof is by induction on order $n$. By Proposition 8, the theorem is trivial if $T$ is a star $S_n$,
hence, in particular, for every tree of order $n=3$.

Suppose that our assertion is true for all trees of order $n-1~(n\geq 4)$ and let T be
a tree of order $n$. We may assume that $T$ is not isomorphic to $S_n$. Let
$x$ be an end vertex of a longest path $P = xyz\dots$ in T and let $T^{'}$ denote the
tree $T-\{x\}$. By the choice of $x$ and $T$, $z$ is the only neighbor of $y$ having
the degree $\geq 2$ in T. Let $d_{T'}(t)$ for any vertex $t\in V(T')$. The degree
in $T'$ of any vertex $t$ in $T'$ is the same as in $T$, except for $t = y$ for which
$d_{T'}(y)=d_{T}(y)-1$.

By induction hypothesis, there is an NFSD-total 2-coloring $f'$ of $T'$. We
will color the edge $xy$ and the vertex $x$ by $a$ and $b$, resp., $a, b\in\{1, 2\}$, so
that the new coloring $f$ of $T$ defined as follows:
$$
f\left( \theta \right) =\left\{ \begin{array}{l}
	f'(\theta)~~~~if~\theta\in V(T')\cup E(T'), \\
	a~~~~~~~~if~\theta=xy,\\
    b~~~~~~~~if~\theta=x,
\end{array} \right.
$$
would be an NFSD-total 2-coloring of $T$. We prove that this is always possible.
Let $\phi'(v)$ denote the expanded sum at $v \in V (T')$ with respect to the coloring $f'$.

Suppose now that the degree $d_{T}(y)$ of $y$ in $T$ is at least 3 and observe
that for any total 2-coloring $f$ of $T$ and for any $t\in N_{T}(y)-\{z\}$ we have
$\phi(t) = f(t)+f(y)+f(yt)\leq 6$ and $\phi(y)\geq 7$, so the vertices $t$ and $y$ are distinguished.
Therefore, we can choose $a$ and $b$ such that $\phi(z) = \phi'(z) \neq \phi'(y)+a +b = \phi(y)$ and
the new total coloring $f$ of $T$ will distinguish all vertices of $T$.

If $d_{T}(y)= 2$, we can also choose $a$ and $b$ such that $\phi(x) = a + f'(y)\neq \phi'(y) + a + b = \phi(y)$
and $\phi(z) = \phi'(z)\neq \phi'(y) + a + b = \phi(y)$, so the total coloring $f$ distinguishes all adjacent vertices of $T$. $\Box$

\vskip 2mm

\noindent {\textbf{Theorem 13.}} \emph{For any hypercube} $Q_n$ $(n\geq 3)$, fgndi$_{\sum}(Q_n)=2$.

\textbf{\emph{Proof}}
Observe that hypercube $Q_n$ is a $n$-regular graph. By Lemma 3, fgndi$_{\sum}(Q_n)\geq 2$.
Next we will prove that $Q_n$ has an NFSD-total $2$-coloring.

For $n=3$ let $$V(Q_3)=\{u_i^{(1)}, v_i^{(1)}| i\in [1,4]\},$$
$$E(Q_3)=\{u_i^{(1)}u_{i+1}^{(1)}, v_i^{(1)}v_{i+1}^{(1)}, u_i^{(1)}v_i^{(1)}| i\in [1,4]\},$$
where subscripts are taken modulo 4.

Let $f_3$ be a total coloring of $Q_3$, and $f_3$ is defined as follows:
$$f_3(u_1^{(1)}u_2^{(1)})=2,f_3(u_3^{(1)}u_4^{(1)})=2,$$
$$f_3(u_2^{(1)}v_2^{(1)})=2,f_3(u_4^{(1)}v_4^{(1)})=2,$$
$$f(\alpha)=1, \alpha\in Q_3-\{u_1^{(1)}u_2^{(1)},u_3^{(1)}u_4^{(1)},u_2^{(1)}v_2^{(1)},u_4^{(1)}v_4^{(1)}\}.$$
Then $f_3$ is an NFSD-total 2-coloring of $Q_3$.

Suppose that $\overline{Q_{n-1}}$ is a copy of $Q_{n-1}$.
Let $$V(Q_{n-1})=\{u_1^{(k)}, u_2^{(k)},u_3^{(k)}, u_4^{(k)}, v_1^{(k)}, v_2^{(k)},v_3^{(k)}, v_4^{(k)}|~k\in [1, 2^{n-4}]\},$$
$$V(\overline{Q_{n-1}})=\{\overline{u_1^{(k)}}, \overline{u_2^{(k)}}, \overline{u_3^{(k)}}, \overline{u_4^{(k)}}, \overline{v_1^{(k)}}, \overline{v_2^{(k)}}, \overline{v_3^{(k)}}, \overline{v_4^{(k)}}|~k\in [1, 2^{n-4}]\}.$$

Observe that each $Q_{n-1}$ and $\overline{Q_{n-1}}$ include
$2^{n-4}$ numbers of $Q_3$, respectively. Meanwhile, $Q_n$ is
constructed in the following procedure, $Q_n=Q_{n-1}\cup
\overline{Q_{n-1}}\cup \{u_i^{k}\overline{u_i^{k}},
v_i^{k}\overline{v_i^{k}}| i\in [1, 4],k\in [1, 2^{n-4}]\}$.

Let $f_n$ be a total coloring of $Q_n$ and be defined recursively as follows.
$$f_n(Q_{n-1})=f_{n-1}(Q_{n-1}),$$ $$f(e_i)=1, i\in [1,2^{n-1}].$$
The total coloring of $\overline{Q_{n-1}}$ is obtained by
exchanging the color 2 of each $Q_3$ in $Q_{n-1}$ (see Fig.3) and
use 1 to color the remaining vertices and edges. $\Box$

\vskip 2mm

\begin{figure}[h]
\centering
\includegraphics[height=3.5cm]{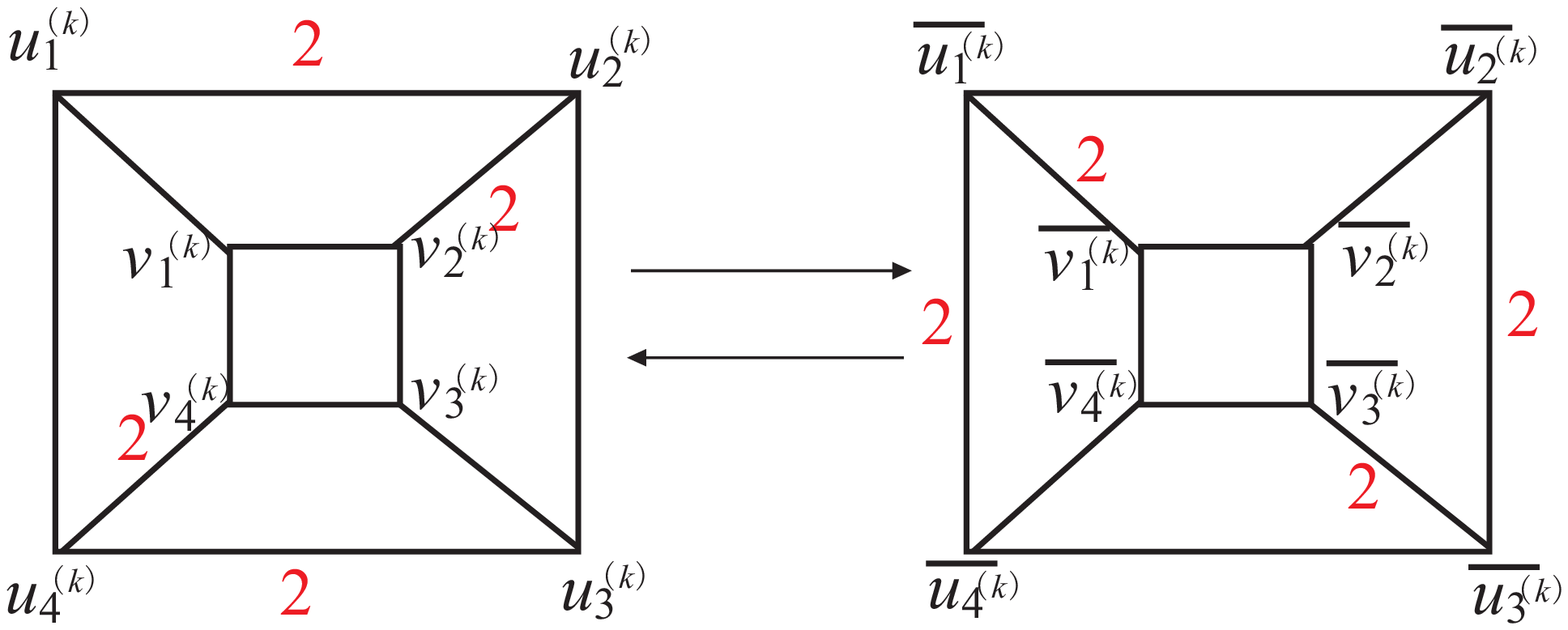}
\begin{center}
{\small Fig.3: Diagram of the exchanging color of a cube.}
\end{center}
\end{figure}

\noindent {\textbf{Theorem 14.}} \emph{Let $G=(X, Y, E)$ be a bipartite graph with
bipartition classes $X$ and $Y$. Then $fgndi_{\sum}(G)\leq 3$}.

\textbf{\emph{Proof}} Let $G=(X, Y, E)$ be a connected bipartite graph with bipartition classes $X$ and $Y$.
If $G$ is a star, then by Proposition 3, the conclusion holds. Now we consider the case that $|X|\geq 2$ and $|Y|\geq 2$.
We define a non-proper total coloring $f$ of $G$ with the following properties:
(1) $f(x)=1$ for any vertex $x$ of $X$;
(2) $f(y)=2$ for any vertex $y$ of $Y$;
(3) Edges between $X$ and $Y$ are colored by $1$.
For an edge $xy$ and $x\in X, y\in Y$, it deduces that $\phi(x)=3d_G(x)+1$ and $\phi(y)=2d_G(y)+2$.
Hence there may appear a case that $\phi(x)=\phi(y)=6k+4$ if $d_G(x)=2k+1$ and $d_G(y)=3k+1$, where $k$ is a positive integer.
Let $N(y)=\{x,v_{x_{1}},v_{x_{2}},\dots,v_{x_{3k}}\}$. Then two cases appear as follows:

\textbf{Case 1.} Not all vertices in $N(y)$ have the same weight $6k+4$. Let $v_{x_{p}}$ be the vertex whose weight is not equal to $6k+4$.

If $\phi(v_{x_{p}})$ is even, recolor $yv_{x_{p}}$ with 2, then the weights of vertices of $G$ keep unchanged as before except for vertices $y$ and $v_{x_{p}}$.
Let $\phi'(y)$ and $\phi'(v_{x_{p}})$ be the new weights of $y$ and $v_{x_{p}}$, respectively. Then $\phi'(y)=\phi(y)+1=6k+5$ and $\phi'(v_{x_{p}})=\phi(v_{x_{p}})+1=3d_G(v_{x_p})+2$.
Thus the weights of $y$ and $v_{x_p}$ are natural distinct. It is easy to verify that $\phi'(v_{x_{p}})=3d_G(v_{x_p})+2\neq 2d_G(v_{y_0})+2=\phi(y_0)$, $y_0\in N(v_{x_p})\setminus \{y\}$. By our coloring rule, for two distinct vertices $v_{x_i}$ and $v_{x_j}$, $|\phi(v_{x_i)}-\phi(v_{x_j})|=6$, so $\phi'(y)\neq \phi(v_{x_0}), v_{x_0}\in N(y)\setminus \{v_{x_p}\}$.

If $\phi(v_{x_{p}})$ is odd, recolor $v_{y}v_{x_{2}}$ with 3, then $\phi'(y)=\phi(y)+2=6k+6$ and $\phi'(v_{x_{p}})=\phi(v_{x_{p}})+2=3d_G(v_{x_p})+3$.
Thus the weights of $y$ and $v_{x_p}$ are natural distinct. It is easy to verify that $\phi'(v_{x_{p}})=3d_G(v_{x_p})+3\neq 2d_G(v_{y_0})+2=\phi(y_0)$, $y_0\in N(v_{x_p})\setminus \{y\}$. Similar to the above discussion, $|\phi(v_{x_i)}-\phi(v_{x_j})|=6$ for two distinct vertices $v_{x_i}$ and $v_{x_j}$, so $\phi'(y)\neq \phi(v_{x_0}), v_{x_0}\in N(y)\setminus \{v_{x_p}\}$.

\textbf{Case 2.} All vertices in $N(y)$ have the same weight $6k+4$.

Recolor edges $xy$ and $yv_{x_{p}} (1\leq p\leq 3k)$ with 2 and recolor vertex $y$ with 1. Then the weight of vertices in $N(y)$ keep the same as before and the weight of $y$ is added to $d_G(y)-1$ than before. Therefore, the weights between vertex $y$ and its neighbors are distinguished. $\Box$

\vskip 2mm

From Theorem 14, the following two results are obvious.

\vskip 2mm

\noindent {\textbf{Corollary 15.}} \emph{Let $G=(X, Y, E)$ be a bipartite graph with
bipartition classes $X$ and $Y$ such that the degree of all vertices of $X$ are even. Then} fgndi$_{\sum}(G)\leq 2$.

\vskip 2mm

\noindent {\textbf{Corollary 16.}} \emph{Let $G$ be a bipartite graph with $\Delta=3$. Then} fgndi$_{\sum}(G)\leq 2$.

 \section{Future Works\ \ }

\vskip 2mm

\noindent {\textbf{Problem 1.}}\ \ \emph{Whether fgndi$_{\Sigma}(G)\leq 3$ holds for every connected graph $G$ with $\Delta=3$} ?

\vskip 2mm

\noindent {\textbf{Problem 2.}}\ \ \emph{Let $K_{n_1,n_2,\dots,n_r}$ be a complete $r$-partite graph with $r$ vertex sets $X_i$ ($i\in [1,r]$) and $|X_i|=n_i$, $\sum_{i=1}^{r}=n$. Besides (i) and (ii) in Theorem 8}, fgndi$_{\Sigma}(K_{n_1,n_2,\dots,n_r})\leq 3$ ?

\vskip 4mm

\bibliographystyle{plain}

\bibliography{ref}  


\noindent {[1]}  L. Addario-Berry, K. Dalal, C. McDiarmid, B.A.
Reed, and A. Thomason, Vertex-colouring edge-weightings[J].
Combinatorica, 2007, 27: 1-12.

\noindent {[2]} L. Addario-Berry, K. Dalal, and B.A. Reed, Degree
constrained subgraphs[J]. Discrete Appl. Math. 2008, 156: 1168-1174.

\noindent {[3]} J.A. Bondy and U.S.R. Murty, Graph theory with
applications. The MaCmillan Press ltd, London and Basingstoke, New
York, 1976.

\noindent {[4]} E. Flandrin, H. Li, A. Marczyk et al. A note on neighbor expanded sum distinguishing index. Discuss. Math. Graph T., 2017, 37(1): 29-37.

\noindent {[5]}  M. Kalkowski, A note on 1,2-Conjecture[D]. Adam Mickiewicz University, 2010.

\noindent {[6]}  M. Kalkowski, M. Karo\'{n}ski, F. Pfender,
Vertex-coloring edge weightings: Towards the 1-2-3-conjecture[J].
J. Combin. Theory Ser. B, 2010, 100: 347-349.

\noindent {[7]} M. Karo\'{n}ski, T. {\L}uczak, A. Thomason, Edge
weights and vertex colours[J]. J. Combin. Theory ser. B, 2004, 91: 151-157

\noindent {[8]} J. Przybylo, M. Wo\'{z}niak, On a 1,2
Conjecture[J]. Discrete Math. Theor. Comput. Sci., 2010, 12(1): 101-108.

\noindent {[9]} J. Przybylo, The 1-2-3 Conjecture almost holds for
regular graphs[J]. J. Combin. Theory ser. B, 2021, 147: 183-200.

\noindent {[10]} T. Wang and Q. Yu, On vertex-coloring
13-edge-weighting[J]. Front. Math. China, 2008, 3: 581-587.

\vskip 1mm

\end{document}